\documentstyle{amsppt}
\magnification1100
\pagewidth{13.5cm}

\define \I{\text {\sl Id}}
\define \Z{\Bbb Z}
\define \R{\Bbb R}
\define \N{\Bbb N}

\define \f{\frac}
\define \bpb{\bigpagebreak}
\define \mpb{\medpagebreak}
\define \spb{\smallpagebreak}
\define \hw{Hantzsche-$\!$Wendt }

\define \ovr{\overset \to}
\define \undr{\underset \to}

\TagsOnRight

\rightheadtext{FLAT MANIFOLDS ISOSPECTRAL ON $\text{p}$-FORMS}
\leftheadtext{R.J.\ Miatello and J.P.\ Rossetti}
\topmatter
\title  FLAT MANIFOLDS ISOSPECTRAL ON $\text{p}$-FORMS \endtitle
\author R.J.\ Miatello  and J.P.\ Rossetti  \endauthor
\address FaMAF-CIEM, Universidad
Nacional  de C\'ordoba, 5000 C\'ordoba, Argentina\. 
\endaddress
\email miatello\@mate.uncor.edu, rossetti\@mate.uncor.edu. \endemail

\thanks 2000 {\it Math Subject Classifications.} primary: 58J53; secondary: 20H15.\endthanks 
\thanks {Partially supported by Conicet and Secyt U.N.C.}\endthanks

\abstract
We study isospectrality on $p$-forms of compact flat manifolds 
by using the equivariant spectrum of the Hodge-Laplacian on the torus. We give an explicit formula for the multiplicity of eigenvalues and a criterion for isospectrality. 
We construct a variety of new isospectral pairs, some of which are the first such examples in the context of compact Riemannian manifolds.
For instance, we give pairs of  flat manifolds of dimension $n=2p , p \ge 2$, not homeomorphic to each other,  which are 
isospectral on $p$-forms but not on $q$-forms for $q\ne p,\, 0\le q \le n$. 
Also, we give manifolds isospectral on $p$-forms if and only if $p$ is odd,  one of them orientable and the other not, and a pair of $0$-isospectral flat manifolds, one of them K\"ahler, and the other not admitting any K\"ahler structure.
We also construct  pairs, $M, M'$  of dimension $n\ge 6$, which are isospectral on functions and such that  $\beta_p(M) < \beta_p(M')$, for $0<p<n$ and
pairs isospectral on $p$-forms for every  $p$ odd, and having different holonomy groups,
$\Z_4$ and $\Z_2^2$ respectively.

\endabstract
\endtopmatter

\document

\heading{\S 1. Introduction}\endheading

This article is a sequel to \cite{MR1} and \cite{MR2}, where  
we use  Sunada's method (\cite{Su}) to produce many pairs of isospectral, non-homeomorphic Hantzsche-Wendt manifolds. 
These are natural generalizations to any odd dimension $n,$ of the classical $3$-dimensional \hw manifold (see \cite{Wo}).

The isospectral manifolds obtained by 
Sunada's method are  always strongly isospectral, hence $p$-isospectral (that is, isospectral on $p$-forms) for all $p$, $0\le p \le n$. The main purpose of the present
article  is to exhibit many examples of compact flat manifolds which are $p$-isospectral for some (but not all) values of
$p$. These examples seem to be new in the context of compact Riemannian manifolds, to our best knowledge.  We will study
$p$-isospectral manifolds by using the equivariant spectrum on the torus, giving an explicit formula for the
multiplicities of the eigenvalues of the Hodge-Laplacian and, as a consequence, a condition for isospectrality on
$p$-forms for each $p$ (see Theorem 3.1).  This formula will also be used to prove non-isospectrality for some flat
manifolds, by computing the multiplicities of specific eigenvalues (see Examples 4.1, 5.1, 5.6). All isospectral
manifolds to be constructed, except those in Example 4.2, will be pairwise non-homeomorphic, since they will have non
isomorphic fundamental groups. 

We consider first  Bieberbach groups with holonomy group 
$\Z_2^r$ and diagonal holonomy representation. 
In this case, the formulas for the multiplicities of eigenvalues involve combinatorial coefficients, namely, integral values of Krawtchouk polynomials,  which do vanish in some cases. 
This allows to produce examples of $p$-isospectral manifolds of dimension $2p$ ($p\ge 2$) which are not $q$-isospectral 
for $q\ne p$. 
We also construct  pairs of manifolds isospectral on $p$-forms for every $p$ odd and having different holonomy
groups: $\Z_4$ and $\Z_2^2$ respectively (see Examples 4.2 and 5.8). We give two 4-dimensional manifolds, $p$-isospectral
for all
$p$, with holonomy group $\Z_2^2$, and having different first integral homology group (see Example 4.5).

In Section 5 we look at pairs  where the holonomy group is $\Z_4 \times \Z_2$ and the holonomy representation is 
not diagonal, giving examples with various properties, in particular:
\spb

-- Two manifolds of dimension $n=6$ which are $0$-isospectral but not isospectral on $p$-forms for any $p,\, 0< p <6$ (Example
5.1; see also Examples 5.4 and 5.5). Moreover one of these manifolds is K\"ahler and the other does not admit any
K\"ahler structure.  By a variation one obtains an isospectral pair where one of the manifolds is hyperk\"ahler and the
other is not (Example 5.3). 

-- Two manifolds which are isospectral on $p$-forms if and only if $p$ is odd, one of them orientable and 
the other not orientable (Example 5.6; see also Example 4.3).

-- Two  $n$-manifolds $M , M'$ which are isospectral on functions but not on $p$-forms for $0<p<n$ and such that $\beta_p(M)< \beta_p(M')$ for each $p,\, 0< p< n$ (Example
5.9).

\mpb

Manifolds which are 0-isospectral and not  $p$-isospectral for some $p$ are not 
very common (see \cite{Pe}). Examples have been given by Ikeda (see \cite{Ik}) for lens spaces, by
Gordon (the first example), Gornet (see \cite{Go},\cite{Gt} respectively) in the context of
nilmanifolds and by Schueth (see \cite{Sch}) for simply connected manifolds.  Flat manifolds yield a very
rich family of non-strongly isospectral pairs with a simple construction and having the property that
certain of their topological invariants can be easily computed in many cases. 
Those constructed in this article are quotients of an 
$n$-torus by free actions of $\Z_2^r$ or $\Z_4\times \Z_2$ and they often yield isospectral manifolds
which are topologically very different from each other, for instance they can be distinguished by their
real cohomology (see Remark 5.10).  The approach in this article can be applied to study isospectrality with
respect to more general differential operators on vector bundles over flat manifolds (see Remark 3.8).

\heading{\S 2. Preliminaries}\endheading

We begin this section by  recalling some known facts on 
Bieberbach groups (see \cite{Ch} or \cite {Wo}).

A discrete, cocompact, torsion-free subgroup $\Gamma$ of $I(\R^n)$ 
is called a Bieberbach group. Then $\Gamma$ acts freely on $\R^n$ and $M_\Gamma = \Gamma\backslash\R^n$ is a compact flat
Riemannian manifold with fundamental group $\Gamma$. 
Any element $\gamma \in
I(\R^n)$ decomposes uniquely $\gamma = B L_b$, with $B \in O(n)$ and $
b\in \R^n.$ 
The  translations in $\Gamma$
form a normal, maximal abelian subgroup  of finite index, $L_\Lambda$,  with $\Lambda$  a lattice in $\R^n$ which is $B$-stable for each $BL_b \in \Gamma$. 
The restriction to $\Gamma$ of the canonical projection
from $I(\R^n)$ to $O(n)$, mapping $BL_b$ to $B$, has kernel $\Lambda$ and the image is a finite subgroup of $O(n)$, called the point group of $\Gamma$. We shall often make the identification of $F := \Lambda\backslash\Gamma$ with the point group of $\Gamma$. The group $F$ coincides with the linear holonomy group of the Riemannian manifold $M_\Gamma$ and the action of $F$ on $\Lambda$ defines an integral representation of $F$, usually called the  holonomy representation. 

We now give a proposition which is useful in the construction of Bieberbach groups. It gives  necessary and
sufficient conditions for a crystallographic group to be torsion-free, in the case of abelian holonomy. 
It is a natural extension of  Proposition 3.1 in \cite{DM}  (valid for $F\simeq {\Z}_2^r$)  and will be used in sections
3 and 4 in the construction of isospectral manifolds.

\bpb

\proclaim{Proposition 2.1}
Assume that  $\Gamma =\langle \gamma_1,\dots,\gamma_r,L_\lambda\,:\,\lambda\,\in\,\Lambda\rangle$  is a subgroup of $Af\!f(\R^n)$, with $\gamma_i = B_i L_{b_i}$, $b_i\in\R^n$, $B_i \in Gl(n,{\R})$ of order $m_i$, 
$\Lambda$  a lattice in $\R^n$ stable by the $B_i$'s  and $\langle B_1,\dots ,B_r\rangle \simeq \Z_{m_1}\times \dots \times \Z_{m_r}$.
We then have that
$\Gamma$ is  torsion-free  with translation lattice
$\Lambda$ if and only if the following two conditions hold:
\roster
\item"(i)" For each pair $i,j,\, 1\le i,j \le r$,  $\,\,\,(B_i^{-1} -\I)b_j - (B_j^{-1}-\I)b_i \in \Lambda$.
\item"(ii)" For  each $I= (i_1,\dots ,i_s)$ 
with $1\le i_1\le i_2\le \dots\le i_s\le r$, let
$B_{i_1}L_{b_{i_1}} \dots B_{i_s} L_{b_{i_s}} = B_I L_{b(I)} \in \Gamma$ ,
with  $B_I:=B_{i_1}\dots
B_{i_s}$ and $b(I)= B_{i_s}^{-1}\dots B_{i_2}^{-1}b_{i_1} + B_{i_s}^{-1}\dots
B_{i_3}^{-1}b_{i_2} + \cdots + B_{i_s}^{-1}b_{i_{s-1}} + b_{i_s}. $
If $B_I$ has order $m$, then 
$$\left(\sum_{j=0}^{m-1} B_I^{-j}\right)b(I) \in \Lambda\,
\setminus \,\left(\sum_{j=0}^{m-1} B_I^{-j}\right)\Lambda.\tag2.1
$$
\endroster

Finally, if $\Gamma$ satisfies conditions (i) and (ii), then   $\Gamma$ is isomorphic to a Bieberbach group with holonomy group $F\simeq \Z_{m_1}\times \dots \times\Z_{m_r}$.
\endproclaim
\demo{Proof}
We first show that conditions (i) and (ii) are necessary for $\Gamma$ to be torsion-free with translation lattice $\Lambda$.
One computes that $[\gamma_i, \gamma_j]=L_{\mu_{i,j}}$,  where
$\mu_{i,j}=B_iB_j (B_j^{-1}-\I){b_i} -(B_i^{-1}-\I){b_j}.$
Since the translation lattice of $\Gamma$ is assumed to be $\Lambda$, this implies that  condition (i) must hold.

Let now $\gamma = B_IL_{b(I)+\lambda}$, with $B_I$, $b(I)$ as in (ii) and $B_I$ of order $m$. Then,
$\gamma^h = B_I^h L_{v_h(I,\lambda)}$, where
$$
v_h(I,\lambda)=
\left(B_I^{-(h-1)}+B_I^{-(h-2)}+\dots +
B_I^{-1}+ \I\right)(b(I)+\lambda)\tag{2.2}
$$ 
Now, since $B_I$ has order $m$, it follows that
$v_m(I,\lambda)$ lies in the translation lattice of $\Gamma$ which is $\Lambda$, and furthermore 
$v_m(I,\lambda)\ne 0$ because $\Gamma$ is torsion-free, hence
condition (ii) must hold. 

Conversely, assume now that (i) and (ii) hold. 
We note that by (i), $\gamma_i \gamma_j =\gamma_j \gamma_i L_{\lambda_{i,j}}$ with $\lambda_{i,j} \in \Lambda$ and also, if $\lambda \in \Lambda, \gamma_i L_\lambda \gamma_i^{-1} =L_{\lambda'}$, for some $\lambda'\in \Lambda$. Furthermore we have that 
$\gamma_i^{-1}= \gamma_i^{m_i -1} L_{\lambda'_i}$ with
$\lambda'_i={-B_i \left(\sum_{j=0}^{m_i -1} B_i^{-j}\right)b_i}\in \Lambda$, by (ii).
Using these facts, one can show that any product of generators in 
$\Gamma$ can be reordered in such a way that a general element  
$\gamma \in \Gamma$ can be written 
$\gamma= \gamma_{i_1}\dots \gamma_{i_s} L_\lambda = B_IL_{b(I)+\lambda} $, for some $I$ as above and
$\lambda \in \Lambda.$ 
Now suppose that $\gamma^h=1$. Then $B_I^h=\I$ and $v_h(I,\lambda) =
0$. This implies, if $m$ is the order of $B_I$, that $h=mk$. Hence 
$\gamma^h = ({\gamma^m})^k= L_{v_m(I,\lambda)}^k = L_{kv_m(I,\lambda)}$,
thus $kv_{m}(I,\lambda)=v_{km}(I,\lambda)=0$ and therefore $v_m(I,\lambda)=0$. 
By (2.2), it now follows that $\left(\sum_{j=0}^{m-1} B_I^{-j}\right)b(I) \in 
\left(\sum_{j=0}^{m-1} B_I^{-j}\right)\Lambda,$
contradicting (ii). This shows that $\Gamma$ is torsion-free.

It only remains to show that the translation lattice of $\Gamma$ equals $\Lambda.$  Any element of $\Gamma$ has the form 
$\gamma= \gamma_1^{l_1}\dots \gamma_r^{l_r} L_\lambda$, with $\lambda \in \Lambda$ and $l_i \ge 0$. If $\gamma =L_\mu$, with $\mu \in \R^n$, then necessarily $l_j =m_j k_j$ for $\, 1\le j \le r$, by the condition on the $B_i$'s in the statement. Now, for each $j$,  $\gamma_j^{l_j}=(\gamma_j^{m_j})^{k_j}=L_{ \lambda_j}^{k_j}= L_{k_j \lambda_j}$, with $\lambda_j \in \Lambda$, by (ii). It thus follows that $\mu \in \Lambda$, as was to be shown. 

Concerning the last assertion, there exists an inner product in $\R^n$ which is invariant under  the holonomy representation of $F$. 
Thus, conjugation by the positive definite transformation relating this inner product to the canonical one,  takes $\Gamma$ into a torsion-free subgroup of $I(\R^n)$. This concludes the proof. 
\qed
\enddemo 

\bpb

\noindent{\bf Remark 2.2.} 
(a) The arguments in the proof actually show that under the assumptions in the statement, $\Gamma$ will have translation lattice $\Lambda$ if and only if condition (i) holds and, for each $1\le i \le r$,  $\left(\sum_{j=0}^{m-1} B_i^{-j}\right)b_i \in \Lambda$.  Furthermore, if this is the case, $\Gamma$ will be torsion-free if and only if  $\left(\sum_{j=0}^{m-1} B_I^{-j}\right)b(I) \notin 
\left(\sum_{j=0}^{m-1} B_I^{-j}\right)\Lambda$, for any $I$ as above.

(b) In  \cite{DM}, for $F\simeq \Z_2^{r},$ it was shown  that the condition for $\Gamma$  to be torsion-free is
$(B_{i_1}\dots B_{i_s} + \I)b(i_1,\dots,i_s)\,\in\,\Lambda \,\setminus \,
(B_{i_1}\dots B_{i_s} + \I)\Lambda$
for any $1\le i_1 < i_2 < \dots < i_s\,\le\,r$.
This condition is equivalent to the condition that the class defined by $\Gamma$ in $H^2(\Z_2^{n-1},\Lambda)$
be a special class (see \cite{Ch, Thm 2.1, p\. 79}).

\bpb

\heading{\S 3. Equivariant spectrum of flat manifolds}\endheading

\spb

If $M$ is a compact Riemannian manifold, let $\text{spec}^p(M)$ denote the spectrum of the Hodge-Laplace 
operator acting on smooth $p$-forms on $M$, $0\le p \le n$. For each $p$, $\text{spec}^p(M)$  is a sequence of
non-negative real numbers tending to $\infty$. Two Riemannian manifolds $M$, $M'$ are said to be {\it
$p$-isospectral} if $\text{spec}^p(M)=\text{spec}^p(M')$. Usually, $0$-isospectral manifolds are just called isospectral. 
Also, $M$, $M'$ are said to be {\it strongly-isospectral} if $\text{spec}_D(M)=\text{spec}_D(M'),$  for any
natural elliptic differential operator $D$ on $M, M'$ (see \cite{DG, Def.\ 3.2}).

Our main goal in this section 
will be  to describe the spectrum of the Laplacian on $p$-forms on $M_\Gamma =\Gamma \backslash \R^n$,
$\Gamma$ a Bieberbach group with translation lattice $\Lambda$. 
We shall use on $\Lambda \backslash{\R}^n$ and on $\Gamma \backslash{\R}^n$ the Riemannian metrics induced by the canonical metric on $\R^n$. 

We will first discuss the function case, $p=0$. 
We note that if $v\in \Lambda ^*$, the dual lattice of
$\Lambda$, the function $f_v (x) = e^{2\pi i\,v.x}$ on ${\R}^n$
is   $\Lambda$-invariant  and $-\Delta f_v = 4\pi ^2 \|v\|^2\,f_v$ ($\Delta$ the Laplacian on ${\R}^n$). Furthermore, 
if $v,v'\in \Lambda ^*$, $v\ne v'$, then 
$\langle f_v ,\,f_{v'} \rangle =0$ and by the Stone-Weierstrass theorem, the set $\{ f_v: v\in \Lambda^*\}$ is a complete orthogonal system of
$L^2(\Lambda \backslash{\R}^n).$ Thus, for each $\mu \ge 0$, 
the eigenspace of $-\Delta$ with eigenvalue $4\pi ^2 \mu$ is
given by ${\Cal H}_\mu := 
 \text{span} \{f_v\,:\, v\in \Lambda^*_\mu\},$ where
$\Lambda^*_\mu = \{v\in \Lambda^*: \| v \|^2=\mu\}$. 
The spectrum of $-\Delta$ in $\Lambda \backslash{\R}^n$ 
is thus determined by the cardinality of the sets 
$\Lambda^*_\mu$.

If we now look at $M_\Gamma = \Gamma\backslash\R^n$, we see that  each 
$\gamma \in \Gamma$ preserves ${
\Cal H}_\mu$ since, for $\gamma =BL_b, B\in O(n), b \in \R^n $, 
$f_v(\gamma x) = e^{2\pi i\,v.(Bx+Bb)} =e^{2\pi i\,B^{-1}v.b} f_{B^{-1}v} (x)$. Hence
$$L^2(\Gamma \backslash{\R}^n) \simeq L^2(\Lambda
\backslash{\R}^n)^\Gamma =  \oplus \sum_\mu {\Cal H}_\mu^\Gamma .$$
Thus $\text{spec}(\Gamma \backslash{\R}^n)=
\{(\mu,d_{\mu})\,:\, \mu\ge 0 \text{ and } d_\mu >0\},$
where $d_{\mu}=\text{dim}\,{
\Cal H}_\mu^\Gamma ,$ for each $\mu.$
\smallskip

The spectrum on $p$-forms on $\Gamma \backslash {\Bbb R}^n$ is obtained in an entirely similar way. 
Let $\omega =\sum_J f_J\,{dx}_J$ be a $p$-form on ${\R}^n$, where
${dx}_J={dx}_{j_1}\wedge\dots\wedge {dx}_{j_p}$ for $J=(j_1,\dots,j_p),\,\, 1\le j_1<\dots<j_p\le n$ and $f_J\in C^\infty (
{\R}^n) $. We have that $\omega$ induces a $p$-form on 
$\Lambda\backslash {\R}^n$  if and only if $L_\lambda^*\,  \omega =\omega$ for any $\lambda \in \Lambda$, that is, 
if $f_J$ is translation invariant by $\Lambda$, for each $J$.

Since $\Delta_p \omega = \sum_J\Delta f_J\, {dx}_J$, it follows that, for each $\mu \ge 0$, 
an orthogonal basis of the eigenspace ${\Cal H}_{p,\mu}$ of
$-\Delta_p$ with eigenvalue $4 \pi^2\mu$ is
 $$\{f_v {dx}_J  : v\in \Lambda^*,\,\|v\|^2= \mu , J=(j_1,\dots,j_p), |J|=p\}.$$  Furthermore, a simple calculation shows that
${\|f_v {dx}_J\|}^2= \text{vol}(\Lambda\backslash \R^n)$, for every $v,J$.

As in the case $p=0$, one has that a form $\omega$ pushes down to $\Gamma\backslash {\R}^n$ if and only if $\gamma ^* \omega =\omega$ for each $\gamma \in \Gamma$, hence  
the eigenspace of $-\Delta_p$ on $p$-forms on $\Gamma \backslash{\R}^n$  
with eigenvalue $4\pi^2 \mu$ is
the space of $\Gamma$-invariants ${\Cal H}^\Gamma _{p,\mu}$, 
provided  ${\Cal H}^\Gamma _{p,\mu}\ne 0$. 
We now give an expression for the dimension $d_{p,\mu}$ of ${\Cal H}^\Gamma_{p,\mu}$, for each $p,\mu.$ 

The map $|F|^{-1}\sum_{\gamma \in \Lambda\backslash\Gamma} 
\,{\gamma^*}_{|{\Cal H}_{p,\mu}} $ is an orthogonal 
projection from ${\Cal H}_{p,\mu}$ onto  ${
\Cal H}_{p,\mu}^\Gamma$, hence 
$d_{p,\mu}=|F|^{-1}\!\left(\sum_{\gamma
\in\Lambda\backslash\Gamma} \, \text{tr} \,{\gamma^*}_{|{\Cal H}_{p,\mu}}
\right)\!.$ 
Now, for each $\gamma=B L_b\in \Gamma$, we have:
$$ \gamma^*(f_v\,{dx}_J)= \gamma^*f_v \,\,B^*{dx}_J=e^{2\pi i\, {B^{-1}}v.b}f_{{B^{-1}}v} \sum_{|J'|=p}c_{J,J'}(B)\, {dx}_{J'}$$
Thus  $\langle \gamma^*(f_v\,{dx}_J)\,,\,f_v\,{dx}_J\rangle
= e^{2\pi i\,v.b}c_{J,J}(B)\,\delta_{Bv,v}\,\text{vol}(\Lambda\backslash \R^n).$ 
Hence
$$\text{tr} \,{\gamma^*}_{|{\Cal H}_{p,\mu}}\!\!= {\text{vol}(\Lambda\backslash
\R^n)}^{-1}
\!\sum_{v:\|v\|^2=\mu }\,\sum _{J:|J|=p}
\langle \gamma^* (f_v {dx}_J) ,f_v {dx}_J \rangle
= \!\!\sum_{\Sb v : \|v\|^2=\mu\\Bv=v\endSb} \!\! e^{2\pi i\,v.b}  
\sum_{|J|=p} c_{J,J}(B) .$$
We note that $\sum_{|J|=p} c_{J,J}(B) =\text{tr}\,\tau_p(B)$, where $\tau_p$ is the canonical representation of $O(n)$ on
$\Lambda^p (\R^n)$. We shall write $\text{tr}_p(B):=\text{tr} \,\tau_p(B)$ ($\text{tr}_0(B)=1$).
 Now for each  $B\in F$ and $\mu \ge 0$, we set 
${{e}}_{\mu, B}(\Gamma):=\!\!\sum_{\Sb v\in \Lambda^*:\|v\|^2=\mu\\Bv=v\endSb} \!\!e^{2\pi i\,v.b}.$ 
Since, for $BL_b\in \Gamma$, $b$ is uniquely determined by $B$ mod $\Lambda$,   ${{e}}_{\mu, B}(\Gamma)$ is well defined.
We have thus proved:

\spb

\proclaim{Theorem 3.1}
If $\Gamma$ is a Bieberbach group with holonomy group $F$, for each $\mu\ge 0$ and $0\le p\le n$, the multiplicity of the eigenvalue $4\pi^2\mu$ of $-\Delta_p$ is given by
$$d_{p,\mu}(\Gamma)=|F|^{-1}\sum_{B \in F} \text{tr}_p(B)\, 
{{e}}_{\mu, B}(\Gamma).\tag{3.1}$$

Let $\Gamma$ and $\Gamma'$ be Bieberbach groups with holonomy groups $F$ and $F'$ and translation lattice $\Lambda$. Let $0\le p\le n$. If there is a bijection $\Phi : B \leftrightarrow B'$ from $F$ onto $F'$ such that, for each $\mu,\, B,$
$$ \text{tr}_p(B)\, {{e}}_{\mu, B}(\Gamma)=\text{tr}_p(B')\, {{e}}_{\mu, B'}(\Gamma')$$
then $M_\Gamma$ and $M_{\Gamma'}$ are isospectral on $p$-forms.
\endproclaim  

\spb

In the examples in this paper we shall always use $\Lambda$ equal to the canonical lattice, hence we will have $\Lambda^*=\Lambda$, $\text{vol}(\Lambda\backslash \R^n)=1$, all $f_v {dx}_J$'s will have norm one and the eigenvalues will be of the form $4\pi^2\mu$ with $\mu \in \N_0$. 

\spb

\example{Remark 3.2}
In the notation of Theorem 3.1 we see that if the bijection $\Phi$ preserves tr$_p$, then 0-isospectral implies
$p$-isospectral. In particular this is always the case if $F=F'$ and if one can take $\Phi=\I$, then 0-isospectral
implies  $p$-isospectral for all $p$. On the other hand, we shall see that $p$-isospectral need not imply
0-isospectral.   
\endexample

\example{Remark 3.3}
If $B\in O(n)$ then tr$_p(B)=\det(B)\text{tr}_{n-p}(B)$, hence (3.1) implies the well
known fact that if $M_\Gamma$ is orientable then $d_{p,\mu}=d_{n-p,\mu}$ for all $\mu\in\N_0$
 (see \cite{BGM, p.\ 238}). 
In particular, if $M_\Gamma$ and $M_{\Gamma'}$ are orientable, then they are $p$-isospectral if and only if they are $(n-p)$-isospectral. 
\endexample

\example{Remark 3.4} We note that if $\mu=0$ we have that 
$$H_{p,0}^\Gamma = {\langle {dx}_J : |J|=p\rangle }^\Gamma = 
{\langle {dx}_J : |J|=p\rangle }^F \simeq
{\Lambda^p(\R^n)}^F,$$
hence $\, d_{p,0}=\dim {\Lambda^p(\R^n)}^F$, which equals 
$\beta_p(M_\Gamma)$ by \cite{Hi, Lemma 2.2}. The equality $\, d_{p,0} = \beta_p(M_\Gamma)$
follows also from the Hodge-de Rham theorem (see \cite{BGM, p.\ 238}).
\endexample

\definition{Definition 3.5}
Given a set $I$ of $h$ elements, and a subset $I_o$ of $I$, with $|I_o|=j$, 
define for $l\le h$, $w_{l,j}(h)= \sum_{\Sb L:L\subset I\\ |L|=l\endSb} 
(-1)^{|L\cap I_o|}.$

It is easy to see that $w_{l,j}(h)$ depends only on $h,j,l$ and not on $I$ and $I_o$. 
Furthermore, for each fixed $t,\,\, 0\le t\le \min(j,l)$, the number of subsets $L \subset I$ such that $|L\cap I_o|=t$
equals $\binom jt \binom{h-j}{l-t}.$ Therefore, we have the formula
$$w_{l,j}(h)= \sum_{t=0}^{\text{min}(j,l)} 
(-1)^t \binom jt \binom{h-j}{l-t}=K_l^h(j).\tag{3.2}$$
where $K_l^h(x):=\sum_{t=0}^{l} 
(-1)^t \binom xt \binom{h-x}{l-t}$ is the  (binary) Krawtchouk polynomial  ($P_l^h(x)$ in the notation of \cite{KL}).

\enddefinition

\example{Remark 3.6}
In the case when the holonomy representation 
diagonalizes in an orthonormal basis of the lattice 
(and we may assume this basis is the canonical basis), 
we can rewrite expression (3.1) in a more explicit form. 
In this case, $Be_i =\pm e_i$, hence $B^* {dx}_J = (-1)^{|I_B' \cap J|}{dx}_J$ where $I_B :=\{i: Be_i =e_i\}$. 
Thus tr$_p(B)= \sum_{J:|J|=p}(-1)^{|J \cap I_B'|}= K_p^n (n-n_B)$, where $I_B'=\{i: i\notin I_B\}$,  $n_B=|I_B|$.
Therefore we have
$$d_{p,\mu}=|F|^{-1}\sum_{v:\|v\|^2=\mu}\, \sum_{\Sb \gamma
\in\Lambda\backslash\Gamma\\Bv=v\endSb} \, K_p^n(n-n_B)\,
e^{2\pi i\,v.b}.  
\tag{3.3} $$
\endexample

The following lemma will be useful. These facts can be found in \cite{KL} or \cite{ChS}, but we include a short proof for completeness. 

\proclaim{Lemma 3.7} If $1\le l,j \le n,$ then we have that
$K_l^n(j) = (-1)^j  K_{n-l}^n(j)$ and \break
$K_l^n(j) = (-1)^l  K_l^n(n-j)$. 
Hence $K_l^n(\frac n2)= K_{\frac n2}^n(j)=0$, for $n$ even and $l,\,j$ odd.
\endproclaim
\demo{Proof} 
Let $I_o \subset I$ be such that $|I_o|=j, \, |I|=n.$
If  $L'$ denotes the complementary subset of $L\subset I$, $|L|=l$, then clearly    
$(-1)^{|L\cap I_o|}=(-1)^{|I_o|}(-1)^ {|L'\cap I_o|}.$
Since $L \rightarrow L'$ gives a bijection between the subsets of $I$ of cardinality $l$ and those of cardinality $n-l$,
it follows that $$K_l^n(j) =  \sum_{\Sb L:L\subset I\\ |L|=l\endSb} (-1)^{|I_o|} (-1)^{|L'\cap I_o|}=
(-1)^j\!
\sum_{\Sb L':L'\subset I\\ |L'|=n-l\endSb} (-1)^{|L'\cap I_o|}= (-1)^j K_{n-l}^n (j).
$$
Similarly, since $(-1)^{|L\cap I_o|}=(-1)^{|L|}(-1)^ {|L\cap I'_o|},$ we have that
$$K_l^n(j)
=\sum_{\Sb L:L\subset I\\ |L|=l \endSb} 
(-1)^{|L|}(-1)^{|L \cap I'_o|}= (-1)^{l} \sum_{\Sb L:L\subset I\\ |L|=l\endSb} 
(-1)^{|L\cap I'_o|}=
 (-1)^l K_l^n (n-j).\qed $$
\enddemo

\example{Remark 3.8}
We  conclude this section by stating a representation theoretic generalization of Theorem 3.1. 

Let $\Gamma$ be a Bieberbach group, $K=SO(n)$ and $G=I(\R^n)$.
Let  V be a finite dimensional inner product space over ${\R}$ or ${\Bbb C}$ and 
let $\rho \times \tau$  a unitary representation of $\Gamma \times K$ on $V$ such that $\rho$ restricts trivially to the
lattice $\Lambda$ of $\Gamma$. One may take $V$ to be the exterior tensor product of a unitary  $\Lambda\backslash \Gamma$
module with a unitary $K$-module (see \cite{Wa, Ch.2}). 

We set $E_{\rho,\tau}= G\times V/\sim\,\,$ where
$(g,v)\sim(\gamma g k, \rho(\gamma)\tau (k^{-1}) v)$, for any 
$\gamma \in \Gamma, \,k\in K, g \in G, v \in V,$
and let $\pi : E_{\rho,\tau} \rightarrow \Gamma \backslash G/K$ be the canonical projection. 
Then $(E_{\rho,\tau}, \pi)$ defines a smooth unitary vector bundle over $\Gamma \backslash G/K,$ 
endowed with a flat connection  and a connection Laplacian $\Delta_{\rho,\tau}$. 
In the special case when  $\rho=1$ and $\tau=\tau_p$ (the exterior representation of $SO(n)$
on $\Lambda^p({\R}^n)$), we get the  exterior vector bundle $\Lambda^p (M_\Gamma)$ over
$M_\Gamma$ and  $\Delta_{\rho,\tau}$ coincides with the Hodge-Laplacian, since the Ricci
tensor is zero.

If  $d_{\rho\times\tau,\mu}(\Gamma)$ denotes the multiplicity of the eigenvalue $4\pi^2 \mu$ of $\Delta_{\rho,\tau}$,
 an argument similar to that in the proof of Theorem 3.1 shows that
$$d_{\rho\times \tau,\mu}(\Gamma)=|F|^{-1}\sum_{\gamma = BL_b \in \Lambda\backslash\Gamma} 
\text{tr}(\rho(\gamma) \tau(B))\, {{e}}_{\mu, B}(\Gamma).\tag{3.4}$$ 
In particular, if $\rho=1$ and   $\tau = \tau_p$, we get the formula (3.1) for the multiplicity 
of the eigenvalues of the Hodge-Laplacian acting on $p$-forms. 
\endexample

\mpb

\heading{\S 4. $p$-Isospectral flat manifolds with diagonal holonomy}\endheading

This section will be devoted to give examples of isospectral flat manifolds with diagonal holonomy, showing in particular that $p$-isospectrality for some $p>0$ need not imply $q$-isospectrality for $q\ne p$ (Example
4.1).  In Example 4.2 we give pairs of $p$-isospectral flat $n$-manifolds ($n\ge 4$) with holonomy group $\Z_2$,  which
can be seen as analogues of the well known pairs of isospectral, non-isometric tori. In 4.4 we discuss the case of \hw
manifolds, recalling  several results on isospectrality proved in \cite{MR1} and \cite{MR2}.  In Example 4.5 we give
a pair of $4$-manifolds, $p$-isospectral for all $p$, whose construction is elementary and having different first
integral homology groups. 

In all cases we shall use Proposition 2.1 in the construction of Bieberbach groups. Condition (i) will hold automatically
since the holonomy representation is diagonal and  $b_j\in \frac12 \Lambda$ (actually the use of Proposition 2.1 in
\cite{DM} would suffice in this section). 

\example{Example 4.1} We will give, for each $p\ge 2$,  a family of Bieberbach groups of dimension $n=2p$, with holonomy
group ${\Z}_2$,  which, generically, are pairwise $q$-isospectral if and only if $q=p$.   

If  k is odd, $1\le k \le n-1$, we set $C_k = \text{diag}(\underbrace{1,\dots ,1}_k , -1,\dots,-1),$ and define $\Gamma_k = \langle C_kL_{\frac {e_1}2}, L_\Lambda \rangle$. 
Clearly, $(C_k+\I)\frac{e_1}2 =e_1 \in \Lambda \setminus (C_k+\I)\Lambda$, hence Proposition 2.1 applies.
We thus get $p$ Bieberbach groups with holonomy group $\Z_2$ which are pairwise non-isomorphic to each other, since the holonomy representations are pairwise not semiequivalent. 
We note that the contribution of the identity element to the sum in (3.3)   
is the same for all $\Gamma_k$'s. Now the  vanishing of the $K_p^n(j)$'s for $j$ odd, by Lemma 3.7, implies that the
contribution of the second element of $F\simeq \Z_2$ to the sum in (3.3) is zero, hence all $\Gamma_k$'s are pairwise
$p$-isospectral.

Furthermore, it is easy to see that these manifolds are not 0-isospectral, for instance by computing $d_1(\Gamma_k)= d_{0,1}(\Gamma_k)$, the multiplicity of the eigenvalue $4\pi^2$.
Indeed, by (3.3) we have that
$$ d_1 (\Gamma_k) =2^{-1}\bigg(2n + \sum_{\Sb v: \|v\|=1\\ C_k v=v \endSb}
 e^{\pi i v.e_1} \bigg) =2^{-1} (2n-2+2(k-1))=n+k-2.$$ 

Similarly we see that the $\Gamma_k$'s are pairwise not $n$-isospectral, by computing $d_{n,1}(\Gamma_k)$. 
In this case by (3.3): 
$$ d_{n,1} (\Gamma_k) =2^{-1}\bigg(2n + \sum_{\Sb v: \|v\|=1\\ C_k v=v \endSb}
K_n^n (|I_{C_k}|) e^{\pi i v.e_1} \bigg) =2^{-1} (2n +(-1)(-2+2(k-1)))=n-k+2.$$

We  note that generically, the $\Gamma_k$'s are pairwise not $q$-isospectral for $q\ne p$. We shall verify this in the case of the pair $\Gamma_1$ and $\Gamma_3$ by showing  that  $\beta_q (\Gamma_1) 
\ne \beta_q (\Gamma_3)$, for  $1\le q \le n-1$, $q\ne p$. 

If $q$ is even, then we compute that
$$\beta_q (\Gamma_1) =\binom{n-1}{q}, \quad \beta_q (\Gamma_3)=
3\binom {n-3}{q-2} + \binom {n-3}{q}.$$
Since $\binom{n-1}{q}= \binom{n-3}{q} + 2\binom{n-3}{q-1} + \binom{n-3}{q-2},$ it follows that $\beta_q (\Gamma_1) = \beta_q (\Gamma_3)$ if and only if
$\binom{n-3}{q-1} = \binom{n-3}{q-2}$, which occurs if and only if $q=\frac n2$.

\spb

In the case when $q$ is odd we have $\beta_q (\Gamma_1) =\binom{n-1}{q-1}, \quad \beta_q (\Gamma_3)=
3\binom {n-3}{q-1} + \binom {n-3}{q-3}.$
An argument entirely similar to the previous one in the case $q$  even, shows that  equality  occurs only if $q=\frac n2$.
Remark 3.4 now implies that $M_{\Gamma_1}$ and $M_{\Gamma_3}$ are not $q$-isospectral for $ 1\le q \le n-1,\,q\ne p$. 
\endexample

\spb

\example{Example 4.2}
The family in Example 4.1 can be enlarged by considering groups $\Gamma_{k,j} 
:= \langle C_kL_{\frac {e_1+\dots+e_j}2}, L_\Lambda \rangle$, for k odd, $1\le j \le k \le n-1$, $n=2p$. 
The $\Gamma_{k,j}$ with $1\le j\le k$, are Bieberbach groups, which for fixed $k$ are pairwise non isometric to each other,  
since they are pairwise not isospectral. Indeed, by (3.3), 
$$d_1(\Gamma_{k,j})=\frac12(2n+2(k-2j))=n+k-2j.$$
However, by arguing exactly as in Example 4.1, we conclude that, for $p=\frac{n}2$, they are all p-isospectral to each
other. 

We shall next describe, for any $n$ even, a  family of $n$-dimensional Bieberbach groups $\Gamma'_h$, all isomorphic to each other, and such that the corresponding manifolds are isospectral on $q$-forms for each $q$ odd, and they are pairwise non isometric. 

For any  $h$ with  $1\le h \le \frac n2$, we set $C:= \text{diag}(\underbrace{1,\dots ,1}_{\frac n2} ,\underbrace{ -1,\dots,-1}_{\frac n2}),$ and define $\Gamma_h := \langle CL_{\frac {e_1 +\dots  +e_h}2}, L_\Lambda \rangle$. We thus get $\frac{n}2$ Bieberbach groups with holonomy group $\Z_2$ which are all isomorphic to each other. Indeed, if  $T_h$ is the linear transformation of $\R^n$ fixing $e_j$, for $j\ge 2$ and such that $T_h e_1 = \sum_{j=1}^h e_j$, then $T_h$ conjugates $\Gamma_1$ onto~$\Gamma_h$.

Since, by Lemma 3.7, $K_q ^n(\frac n2)=0$ for any $q$ odd,  it follows that the associated flat manifolds are isospectral
on $q$-forms for $q$ odd, by (3.3). One can again see that they are not isospectral on functions (hence they are not
isometric) by computing the multiplicity of the first non zero eigenvalue.  These manifolds are quotients of a flat torus
of dimension $n$ by  free actions of $\Z_2$. In particular, in dimension  $n=4$, we get two very simple manifolds which
are diffeomorphic, non isometric, and $q$-isospectral for $q=1,3$. We note that these examples can be seen as analogues
of the well known isospectral, non isometric tori, which exist for $n\ge 4$ (see \cite{CS} and references therein).
\endexample

\example{Example 4.3}
We now construct a pair of flat 9-manifolds which are $p$-isospectral only for one odd value of $p$ and one even value of $p$. Moreover,  one of them will be orientable and the other not orientable.
We use $n=9$ since it is the first odd dimension such that some of the coefficients $K_l^n(j)$ vanish.
Actually, one computes easily using (3.2) that $K_l^9(j)=0$ for $(l,j)=
(2,3);\, (2,6);\,(3,2);\,(3,7);\,(6,2);\,(6,7);\,(7,3)$ and $(7,6)$.

We let $B:= \text{diag}(1,1,1,-1,-1,-1,-1,-1,-1)$, 
$B':= \text{diag}(1,1,1,1,1,1,-1,-1,-1)$.
We define $\Gamma := \langle BL_{\frac {e_1}2}, L_\Lambda \rangle$ and  $\Gamma' := \langle B'L_{\frac {e_1}2}, L_\Lambda \rangle$, $\Lambda$ the canonical lattice. As in Example
4.1 we see that Proposition 2.1 applies. We note that $M_{\Gamma}={\Gamma}\backslash \R^9$ is orientable but
$M_{\Gamma'}={\Gamma'}\backslash \R^9$ is not orientable.

The vanishing of $K_2^9(6)$ and $K_7^9(6)$ imply that
$\text{tr}_2(B)=\text{tr}_7(B)=0$. Similarly, 
$\text{tr}_2(B')=\text{tr}_7(B')=0$ since $K_2^9(3)=K_7^9(3)=0$.
It follows from Theorem 3.1 that $M_\Gamma$ and $M_{\Gamma'}$ are isospectral on $p$-forms for $p=2\text{ and }7$.

The non isospectrality for the other values of $p$ can be seen as follows.
By Remark 3.4, the computation of the Betti numbers shows the non $p$-isospectrality for $p=1,3,5$ and 9. 
Indeed, one computes that the Betti numbers are: 1, 3, 18, 46,   60,   60,  46,   
18,   3,  1 for $M$ and 1,  6,   18,   38,   60,   66,  46,    18,   3, 0 for $M'$,
respectively.

For the remaining values of $p$ (that is, 0,4,6 and 8) one can argue as in Example 4.1, showing that 
$d_1(\Gamma) < d_1(\Gamma')$, which proves non isospectrality on functions. 
Since $\text{tr}_p(B)=\text{tr}_p(B')\ne 0$ for $p=4,6$ and 8, this implies that $d_{p,1}(\Gamma) \ne 
d_{p,1}(\Gamma')$ for these values of $p$.

\endexample

\example{Example 4.4}
An interesting case of flat manifolds with diagonal holonomy representation is that of \hw manifolds. These are compact flat $n$-manifolds with holonomy group $\Z_2^{n-1}$, which are rational homology spheres.
We will recall some basic facts from \cite{MR1}. Let $n$ be odd. A {\it \hw  group} (or {\it HW   group})  $\Gamma$ is an $n$-dimensional Bieberbach group with holonomy group $F\,\simeq\,\Z_2^{n-1}$ and such that  the action of every $B\in F$ diagonalizes on a $\Z$-basis $v_1,\dots, v_n$  of $\Lambda$ with $\det(x)= 1.$ The holonomy group $F$ can thus be identified to  the diagonal subgroup 
$\{B:Bv_i = \pm v_i, \,\,1\,\le\,i\,\le\,n,\,\,\det B = 1\}$.
Since $n \ge3$, for each pair $j \ne k$, there exists $B \in F$ such that $Bv_j =v_j,\, Bv_k=-v_k$, therefore it follows that all the vectors $v_i$ must be orthogonal to each other. Hence,
after conjugation of $\Gamma$ by an element in $O(n,\R)$ we may change
the translation lattice of $\Gamma$ into a lattice generated by  vectors $c_ie_i$,  $1\le i \le n$, where $e_1,\dots, e_n$ is the canonical basis of ${\R}^n$. For the purposes of constructing isospectral flat manifolds, we will assume as usual, that $\Lambda = \sum_{i=1}^n \Z e_i,$  the canonical lattice, which is the one which allows more symmetries.
The corresponding manifold $M_\Gamma :=\Gamma\backslash\R^n$ is
called a {\it Hantzsche-Wendt} (or {\it HW}) {\it  manifold}.

We denote by $B_i$ the diagonal matrix  fixing $e_i$ and such that  $B_i e_j = -e_j$ (if $j \ne i)$, for each $1\,\le\,i\,\le\,n$.
Clearly, $F$ is generated by $B_1, B_2,\dots, B_{n-1}$. Furthermore 
$B_n = B_1 B_2 \dots B_{n-1}$.

If $\Gamma$ is an HW group, then $\Gamma = \langle B_1
L_{b_1},\dots,B_{n-1}L_{b_{n-1}},L_\lambda\,:\,\lambda\,\in\,\Lambda\rangle$, 
for some $b_i\in\R^n$, $1\le i \le n-1$, where it may be assumed that 
the components of $b_i$, $b_{ij}\in\{0,\frac12\}$, for $1\,\le\,i,j\,\le\,n.$ The matrix $A := [b_{ij}]$  plays a main 
role in the study of HW groups. 
If  $I= \{i_1,\dots, i_s\}$, with $1\le i_1<\dots <i_s \le n$, we write $B_I= B_{i_1}\dots B_{i_s}.$  Since $b_{ij}\,\in\,\{0,\frac12\}$ for all $i,j$, it follows that any $\gamma \in \Gamma$ can be written {\it uniquely} $\gamma= B_I L_{b_o(I)}L_\lambda$
with $s$ odd, $\lambda \in \Lambda$ and $b_o(I)_j \in \{0,\f1{2}\},$ for
$1\le j\le n $. Condition (i) in Proposition 2.1 is automatic, while (ii) says that $(B_I +
\I)b_o(I)\,\in\,\Lambda\setminus (B_I + \I)\Lambda,$ for any $I$ as above.

We now  rewrite  formula (3.3) for the multiplicities of eigenvalues 
$d_{p,\mu}$ in the present case.
For each $v\in \Lambda$ we set
$I_v :=\{ j: v.e_j\ne 0 \},\quad I_v^{odd} :=\{ j:v.e_j \text{ is
odd}\}.$
One has that $B_I e_j =-e_j$ if and only if $j \notin I$, since $|I|$ is odd. 
Hence $B_I v =v$ if and only if  $I \supset I_v$. 
Furthermore, $e^{2\pi i v.b_o(I)} = (-1)^{|I_{2b_o(I)} \cap I_v^{odd}|}$ and, 
by Lemma 3.7, $K_p^n(n-n_B)=(-1)^p K_p^n(n_B)$. 
Therefore, for each $\mu\in \N,\, 0\le p \le n$ we obtain
$$
d_{p,\mu}=
{\tsize\frac{{(-1)}^{p}}{|F|}}
\sum_{\|v\|^2=\mu}\,\,\sum_{\Sb I:I\supset I_v\\ |I| \,\, odd \endSb} 
K_p^n(|I|)(-1)^{|I_{2b_o(I)}\cap I_v^{odd}|}. \tag{4.1}
$$
Formula (4.1) can be used to construct many isospectral pairs of \hw manifolds. 
In \cite{MR1} and \cite{MR2} by using Sunada's theorem,  
we have shown that the number of pairs of isospectral HW  manifolds non homeomorphic to each other grows exponentially with
$n$. In \cite{MR2} we showed that there are exactly $62$ diffeomorphism classes of HW manifolds in dimension $7$, giving a
full set of representatives and using (4.1) to list all isospectral classes for $p=0$. {}From this classification it
follows in particular that isospectrality implies Sunada-isospectrality, hence  $p$-isospectrality for all $p$, for HW
manifolds in dimension
$n=7$. Also, there is evidence that there exist large isospectral sets, for instance there
are several families consisting of 10 pairwise isospectral HW manifolds in dimension 9.

In the more general setting of flat manifolds with diagonal holonomy representation, we have proved by combinatorial methods that
$0$-isospectrality implies $p$-isospectrality and,
generically, $p$-isospectrality for one value of $p$ implies $p$-isospectrality for all $p$. As we have seen, there are
some exceptions to this rule, due to the vanishing of some of the combinatorial numbers $K_l^n(j)$, which allow to have
$p$-isospectrality for certain values of $p$ only. Even in the case of HW manifolds, there exist manifolds of dimension 9
which are isospectral on $3$-forms but not $0$-isospectral.  We discuss these topics in another paper (see \cite{MR3}). 
\endexample

\example{Example 4.5} 
Among the known examples of isospectral, non homeomorphic  manifolds in low dimensions 
one can mention those in \cite {Vi} (hyperbolic 3-manifolds), \cite{Ik} (lens spaces, for $n$ odd, $n\ge 5$), \cite{Go} 
(nilmanifolds, $n\ge 5$)  and \cite{DM} (flat manifolds with $n\ge 5$). We shall now give a pair of isospectral flat manifolds
$M_\Gamma$, $M_{\Gamma'}$ of dimension $n=4$, with holonomy group $\Z_2^2$,  whose construction is elementary. We will give
their Betti numbers, and will show that $H_1(M_\Gamma,\Z)$ and $H_1(M_{\Gamma'},\Z)$ are not isomorphic. 

Now let $B_1=\undr{ \left[ \smallmatrix 1& & & \\ &1& & \\ & &-1& \\ & & &-1 \endsmallmatrix \right]}$, $B_2=\left[
\smallmatrix 1& & & \\ &-1& & \\ & &1& \\ & & &-1 \endsmallmatrix \right]$, $b_1=\frac{e_2+e_4}2$, $b_2=\frac{e_3}2$,
$b_1'=\frac{e_2}2$, $b_2'=\frac{e_1}2$,  and let $\Gamma=\undr{\langle B_1L_{b_1}, B_2L_{b_2}, \Z^4\rangle}$ and  
$\Gamma'=\undr{\langle B_1L_{b_1'}, B_2L_{b_2'}, \Z^4\rangle}$.
One has that $B_1L_{b_1}\, B_2L_{b_2}=B_1B_2L_{\frac{e_2+e_3+e_4}2}\, L_\lambda$, with $\lambda=-e_2-e_4 \in \Lambda$. One easily sees that 
condition (ii) in Proposition 2.1 holds, hence
$\Gamma$ and $\Gamma'$ are  Bieberbach groups (they correspond respectively to those  denoted by 5/1/2/7 and 5/1/2/9 in \cite{BBNWZ}).

A direct computation shows that $H_1(M_\Gamma,\Z)\simeq \Gamma / [\Gamma,\Gamma]= \Z\oplus \Z_4^2$ and $H_1(M_\Gamma',\Z)\simeq \Gamma' / [\Gamma',\Gamma']= \Z\oplus \Z_2^3$ (see alternately the table in \cite{RT} \S 7, the manifolds with parameters $r=m_1=m_2=m_3=1$ and special classes  $(0,1,1,1)$ and $(h_1,1,0,0)$ res\-pec\-ti\-vely). 
The Betti numbers of  $M_\Gamma$ and $M_\Gamma'$ are  immediately computed by using that $\beta_j=\text{dim}{\Lambda^j(\R^4)}^F$ for $0\le j\le 4$ (Remark
3.4).  One has that ${\Lambda(\R^4)}^F=\nobreak{\langle 1, e_1, e_2\wedge e_3\wedge e_4, e_1\wedge e_2\wedge e_3\wedge
e_4 \rangle}$, hence $\beta_j=1$ for $j\neq 2$ and $\beta_2=0$, for both, $M_\Gamma$ and $M_\Gamma'$.

We shall verify next, using  Theorem 3.1 with $\Phi=\I$, that $M_\Gamma$ and $M_{\Gamma'}$ are isospectral. 
We first see that for each $\mu\in\N_0$ and $B\in F$, 
${{e}}_{\mu,B}(\Gamma)={{e}}_{\mu,B}(\Gamma')$. 
By Remark 3.2 it will then follow that $M_\Gamma$ and $M_\Gamma'$ are $p$-isospectral for all $p$.

In the first place, if $B=\I$, one always has that ${{e}}_{\mu,Id}=\text{card}\{v:{||v||}^2=\mu\}$. If $B\neq\I$ one has:
$${{e}}_{\mu,B_1}(\Gamma)=\sum_{\Sb {||v||}^2=\mu \\ v\in\langle e_1, e_2\rangle \endSb} e^{2\pi i \frac{e_2+e_4} 2 \cdot v }=
\sum_{\Sb {||v||}^2=\mu \\ v\in\langle e_1, e_2\rangle \endSb} e^{2\pi i \frac{e_2} 2 \cdot v }={{e}}_{\mu,B_1}(\Gamma');$$
$${{e}}_{\mu,B_2}(\Gamma)=\sum_{\Sb {||v||}^2=\mu \\ v\in\langle e_1, e_3\rangle \endSb} e^{2\pi i \frac{e_3} 2 \cdot v} =
\sum_{\Sb {||v||}^2=\mu \\ v\in\langle e_1, e_3\rangle \endSb} e^{2\pi i \frac{e_1} 2 \cdot v }={{e}}_{\mu,B_2}(\Gamma'),$$
where the central equality of the second line holds since the sum is symmetric with respect to  $e_1, e_3$. 
In the case of $B_1$ we have used the fact that one can always disregard any vector perpendicular to the  space fixed by $B$ (in this case $\frac {e_4}2$). 
The equality in the case of $B_1B_2$ results  from  combining these two observations. Indeed 
$${{e}}_{\mu,B_1B_2}(\Gamma)=\sum_{\Sb {||v||}^2=\mu \\ v\in\langle e_1, e_4\rangle \endSb} e^{2\pi i \frac{e_2+e_3+e_4} 2 \cdot v }=
\sum_{\Sb {||v||}^2=\mu \\ v\in\langle e_1, e_4\rangle \endSb} e^{2\pi i \frac{e_1} 2 \cdot v }={{e}}_{\mu,B_1B_2}(\Gamma').$$

We observe that  one can  show that the manifolds in the present example are actually Sunada-isospectral, 
by using the same methods as in \cite{MR1, \S 3}, for instance.  

\endexample

\head \S 5. Isospectral, non-strongly  isospectral flat manifolds  
\endhead

In the present section we shall construct several examples of pairs of  non homeomorphic manifolds  with holonomy group $\Z_4\times\Z_2$,  which are isospectral on $p$-forms for some $p,\, 0\le p \le n$, but not for all values of $p$. 
The use of non-diagonal holonomy representation will allow to construct various isospectral pairs with new properties, as described in the Introduction. The procedure in most cases will be to define a bijection $\Phi\colon F\rightarrow F'$ as in  Theorem
3.1, which does not preserve the value of tr$_p(B)$, for some values of $p$.

\spb

\example{Example 5.1}
Let $\widetilde J = \left[ \smallmatrix 0 & 1 \\ -1 & 0 \endsmallmatrix \right]$ and let $\I_2$ be the identity transformation  on $\R^2$. One has that 
${\widetilde J}\, ^4 =\I_2$. 

Let $\Gamma = \langle B_1L_{b_1}, B_2L_{b_2}, \Lambda \rangle$, and 
$\Gamma' = \langle B'_1L_{b'_1}, B'_2L_{b'_2}, \Lambda \rangle$, where 
$\Lambda = \Z^6$ and 
$$
B_1=\left[ \matrix \widetilde J & & \\  & \widetilde J & \\ & & \I_2 \endmatrix \right] \qquad
B_2=\left[ \matrix -\I_2 & &  \\  & \I_2 & \\ & &  \I_2 \endmatrix \right]
$$
$$
B'_1=\left[ \smallmatrix \widetilde J & & & & \\  & 1 & & & \\ & & -1 & & \\ & & & -1 & \\ & & & & 1 \endsmallmatrix \right] \qquad
B'_2=\left[ \smallmatrix -\I_2 & & & & \\  & -1 & & & \\ & & 1 & & \\ & & & -1 & \\ & & & & 1 \endsmallmatrix \right],
$$
$$b_1=\frac{e_5}4, \,\, b_2=\frac{e_6}2,\,\,\,\, b'_1=\frac{e_6}4,\,\, b'_2=\frac{e_4+e_5}2.$$ 

We first verify conditions (i) and (ii) in Proposition 2.1 to ensure that $\Gamma, \Gamma'$ are Bieberbach groups.
By the choices of $b_1,b_2$, for any $BL_b ,\, CL_c \in \Gamma$, one has $Cb=b$ and this implies that (i) holds for $\Gamma$.
In the case of $\Gamma'$ we have $({B'_2}^{-1}-\I)b'_1  - ({B'_1}^{-1}-\I)b'_2 = e_4 +e_5 \in \Lambda$.

Condition (ii) holds for $\Gamma$ since $\left(\sum_{j=0}^{m-1} B^{-j}\right)b =mb \in \Lambda \setminus \left(\sum_{j=0}^{m-1} B^{-j}\right)\Lambda$. For instance, if $B=B_1$, we get $\left(\sum_{j=0}^{3} B_1^{-j}\right)b_1=e_5 \in \Z^6\setminus 4\Z e_5 \oplus 4\Z e_6$.
In the case of $\Gamma'$ one can argue in a similar way, using the fact that $B'e_6 = e_6$ for any $B'L_{b'} \in \Gamma'$.

It follows that $\Gamma$ and $\Gamma'$ are Bieberbach groups with translation lattice $\Z^6$ and holonomy group $\Z_4\times\Z_2$.
Since the holonomy representations  are not semiequivalent to each other,  $\Gamma$ and $\Gamma'$ are not isomorphic (see \cite{Ch, p. 81}).

One has that 
$F=\{ B_1^iB_2^j\colon 0\le i\le 3,\, 0\le j\le 1 \},\,\,\, 
F'=\{ {B'_1}^i{B'_2}^j\colon 0\le i\le 3,\, 0\le j\le 1 \}$.
In the following table we give, in a  visual way, a list of the non trivial representantives 
of $\Lambda \backslash\Gamma$ and $\Lambda \backslash\Gamma'$.
One represents each element of $F$ and $F'$ by a column, indicating the  (non zero) translational components modulo $\Lambda$ by means of subindices. For example $B_1L_{b_1}$ will be represented by  the column  ${(\widetilde J, \widetilde J, 1_{\frac 14}, 1)}^t$, where $\frac 14$ in the fifth component indicates that $b_1=\frac {e_5}4$.
As each translational component  is only determined modulo $\Lambda$, one may choose it, in each  case, in such a way that its coordinates lie in $[0,1)$.

\

\centerline{
\vbox{\tabskip=0pt
\def\tablerule{\noalign{\hrule}}
\halign to340pt {\strut#& \vrule#\tabskip=1em plus 2em&
 \hfil# \hfil&  #& 
 \hfil# \hfil&  #& 
 \hfil# \hfil& \vrule#& 
 \hfil# \hfil& \vrule#& 
 \hfil# \hfil&  #&  
 \hfil# \hfil&  #& 
 \hfil# \hfil&  \vrule#  
\tabskip=0pt\cr\tablerule
&& {\vrule height14pt width0pt depth8pt} $B_1$ && $B_1^2$ && $B_1^3$ && $B_2$  && $B_1B_2$ && $B_1^2B_2$ &&  $B_1^3B_2$ & \cr\tablerule 
&& {\vrule height12pt width0pt depth0pt}  $\widetilde J$  && $-\I_2$  && $-\widetilde J$ && $-\I_2$ && $-\widetilde J$ && $\I_2$ && $\widetilde J$ &\cr
&&  $\widetilde J$  && $-\I_2$  && $-\widetilde J$ && $\I_2$ && $\widetilde J$ && $-\I_2$ && $-\widetilde J$ &\cr
{\vrule height0pt width0pt depth8pt}&& $1_{\frac 14}$ && $1_{\frac 12}$ && $1_{\frac 34}$ && 1 && 
$1_{\frac 14}$ && $1_{\frac 12}$ && $1_{\frac 34}$ &\cr
{\vrule height0pt width0pt depth8pt}&& 1 && 1 && 1 && $1_{\frac 12}$ && $1_{\frac 12}$ && $1_{\frac 12}$ && $1_{\frac 12}$ &\cr
\tablerule
\cr}}
}

\

\centerline{
\vbox{\tabskip=0pt
\def\tablerule{\noalign{\hrule}}
\halign to348pt {\strut#& \vrule#\tabskip=1em plus 2em&
 \hfil# \hfil&  #& 
 \hfil# \hfil&  #& 
 \hfil# \hfil& \vrule#& 
 \hfil# \hfil& \vrule#& 
 \hfil# \hfil&  #&  
 \hfil# \hfil&  #& 
 \hfil# \hfil&  \vrule#  
\tabskip=0pt\cr\tablerule
&& {\vrule height14pt width0pt depth8pt} $B'_1$ && ${B'_1}^2$ && ${B'_1}^3$ && $B'_2$  && $B'_1B'_2$ && ${B'_1}^2B'_2$ &&  ${B'_1}^3B'_2$ & \cr\tablerule 
&& {\vrule height12pt width0pt depth0pt}  $\widetilde J$  && $-\I_2$  && $-\widetilde J$ && $-\I_2$ && $-\widetilde J$ && $\I_2$ && $\widetilde J$ &\cr
&& 1 && 1  && 1 && $-1$ && $-1$ && $-1$ && $-1$ &\cr
{\vrule height0pt width0pt depth8pt}&& $-1$ && 1  && $-1$ && $1_{\frac 12}$ && $-1_{\frac 12}$ && $1_{\frac 12}$ && $-1_{\frac 12}$ &\cr
{\vrule height0pt width0pt depth8pt}&& $-1$ && $1$ && $-1$ && 
$-1_{\frac 12} $ && $1_{\frac 12}$ && $-1_{\frac 12}$ && 
$1_{\frac 12}$ &\cr
{\vrule height0pt width0pt depth8pt}&& $1_{\frac 14}$ && $1_{\frac 12}$ && $1_{\frac 34}$ && $1$ && $1_{\frac 14}$ && $1_{\frac 12}$ && $1_{\frac 34}$ &\cr
\tablerule
\cr}}}
\

We now analyze isospectrality by using Theorem 3.1. We define a bijection $\Phi \colon F\rightarrow F'$ by
$\Phi({B_1}^i{B_2}^j)= {B'_1}^i{B'_2}^j$, if $i=1,3$ or $i=2, j=1$ and by $\Phi(B_1^2)=B'_2,\,\, \Phi(B_2)={B'_1}^2$. One
observes that $\Phi$ preserves, except for a permutation of the coordinates, the spaces fixed by the matrices $B_i$  and
the projections of the translational components onto these fixed spaces. For example, since ker$(B_1-\I)=\langle e_5,
e_6\rangle,\, b_1=\frac {e_5}4$ and for $B'_1=\Phi(B_1)$ one has  ker$(B'_1-\I)=\langle e_3, e_6\rangle,\, b'_1=\frac
{e_6}4$. One observes that the spaces fixed by $B_1$ and $B'_1$  are the same up to a permutation map which sends $e_3$
to $e_5$, $e_5$ to $e_6$, $e_6$ to $e_3$ and leaves the other $e_i$'s fixed. In the case of $B_1^2$ and
$B'_2=\Phi(B_1^2)$ one has  ker$(B_1^2-\I)=\langle e_5, e_6\rangle$, and the  translational component of $B_1^2$ is
$\frac {e_5}2$ while   ker$(B'_2-\I)=\langle e_4, e_6\rangle$ and $b'_2=\frac {e_4+e_5}2$ which, when projected onto the
fixed space  of $B'_2$ yields $\frac {e_4}2$. Taking into account these considerations, with an argument analogous to
that of Example 4.5, one verifies that 
${{e}}_{\mu,B}(\Gamma)={{e}}_{\mu,\Phi(B)}(\Gamma'),\,\,\, \forall B\in F$, hence $M_\Gamma$ and $M_{\Gamma'}$ are isospectral on functions.
Since they are orientable, it follows that they are also $6$-isospectral.

We now see that $M_\Gamma$ and $M_{\Gamma'}$  are not $p$-isospectral, for $1\le p\le 5$.
In the first place we  observe that $B_1^2$ and $B'_2$ are conjugate by a permutation matrix, and the 
same happens with $B_2$ and ${B'_1}^2$ and 
with $B_1^2B_2$ and ${B'_1}^2B'_2$. As a consequence, tr$_p(B_1^2)=\text{tr}_p(B'_2)$, $0\le p\le 6$, hence, 
since ${{e}}_{\mu,B_1^2}(\Gamma)={{e}}_{\mu,B'_2}(\Gamma')$,
it follows that $B_1^2$ and $B'_2$ give the same contribution to the sum in (3.1), for each $p$. The same is true for
the  two remaining pairs. 

Since $\widetilde J$ and $-\widetilde J$ are conjugate it turns out that all 4 elements of order 4 
in $F$ (resp\. $F'$): $B_1, B_1^3, B_1B_2$ and $B_1^3B_2$ (resp\. $B'_1, {B'_1}^3, B'_1B'_2$ and ${B'_1}^3B'_2$) are
conjugate to each other, hence their $p$-traces  are the same, for each $p$. The $p$-traces of $B_1$ and $B'_1$ are
respectively given by 

\

\centerline{
\vbox{\tabskip=0pt
\def\tablerule{\noalign{\hrule}}
\halign to240pt {\strut#& \vrule#\tabskip=1em plus 2em&
 \hfil# \hfil&  \vrule#& 
 \hfil# \hfil&  #& 
 \hfil# \hfil&  #& 
 \hfil# \hfil&  #& 
 \hfil# \hfil&  #& 
 \hfil# \hfil&  \vrule#  
\tabskip=0pt\cr\tablerule
&&  $p$\kern -0pt && $1$ && $2$ && $3$  && $4$ && $5$ & \cr\tablerule 
&&  tr$_p(B_1)$ \kern -6pt && $2$  && $3$ && $4$ && $3$ && $2$  &\cr\tablerule
&& tr$_p(B'_1)$ \kern -6pt&& $0$ && $-1$ && $0$ && $-1$ && $0$ &\cr
\tablerule
\cr}}
}

\

\noindent For instance, if $p=2$, $B_1$ induces a linear transformation in $\Lambda^2(\R^6)$ such that the associated matrix in the basis $\{ e_i\wedge e_j : i<j \}$ has only three non zero diagonal entries corresponding to $\,e_1\wedge e_2,\, e_3\wedge e_4\,$ and $\,e_5\wedge e_6$, all of them equal to 1. Hence tr$_2(B_1)=3$. However, in the case of $B'_1$ the non zero diagonal entries correspond to $\,e_1\wedge e_2$,  $e_3\wedge e_6$, $e_4\wedge e_5$ (equal to $1$) and to $\,e_3\wedge e_4$,  $e_3\wedge e_5$, $e_4\wedge e_6\,$ and $\,e_5\wedge e_6$ (equal to~$-1$). Hence, tr$_2(B'_1)=-1$. The remaining $p$-traces for $B_1$ and $B'_1$ are computed similarly. Furthermore, by Remark
3.3,\  tr$_p(B)=\det (B) \,\text{tr}_{n-p}(B)$ for any $\,B\in \text{O}(n)$, hence in general it suffices to compare the
multiplicities $d_{p,\mu}$ for $\,p\le \left[\frac n2\right]$.

By (3.1) the expressions of $d_{p,\mu}(\Gamma)\,$ and $\,d_{p,\mu}(\Gamma')\,$ are as follows:
$$\align
\!d_{p,\mu}(\Gamma)\!&=\!|F|^{-1}\!\!\left(\text{tr}_p(\I) {{e}}_{\mu,Id}(\Gamma) + \!\!\sum_{B\colon B\neq Id} \!\text{tr}_p(B)\, {{e}}_{\mu,B}(\Gamma)\right)= \\
&=\!|F|^{-1}\!\!\left(\!\binom np \left|\{v\colon {\|v\|}^2=\mu\}\right|\! + 
\!\text{tr}_p(B_1)\!\!\!\! \sum_{B\colon \text{ord}(B)=4} \!\!{{e}}_{\mu,B}(\Gamma)\! + \!\!\!\!
\sum_{B\colon \text{ord}(B)=2} \!\!\!\!\text{tr}_p(B) {{e}}_{\mu,B}(\Gamma)\!\!\right) \\
d_{p,\mu}(\Gamma')\!&=\!|F'|^{-1}\!\!\left(\!\!\!\binom np \!\left|\{v\colon {\|v\|}^2\!=\!\mu\}\right|\! +\! 
\text{tr}_p(B'_1) \!\!\!\!\!\!\sum_{B'\colon \text{ord}(B')=4} \!\!\!\!\!{{e}}_{\mu,B'}(\Gamma')\! +\!\!\!\!\!\!\!
\sum_{B'\colon \text{ord}(B')=2}\!\!\!\!\!\!\! \text{tr}_p(B') {{e}}_{\mu,B'}(\Gamma')\!\!\right)\!\!. 
\endalign
$$
By  the previous observations it turns out that the only possible difference between both expressions is in the second terms. Since   
$\sum_{\text{ord}(B)=4} {{e}}_{\mu,B}(\Gamma) =
\sum_{\text{ord}(B)=4} {{e}}_{\mu,B'}(\Gamma')\,$ and, as seen above,
tr$_p(B_1)\neq \text{tr}_p(B'_1)$, for $\,1\le p\le 5$, the proof will be complete if we verify
that $\sum_{\text{ord}(B)=4} {{e}}_{\mu,B}(\Gamma)\neq 0$, for some value of $\mu$. The  fixed space of each $\,B\in F$ with ord$(B)=4\,$ is $\,\langle e_5, e_6 \rangle$. We take $\mu=8$. The vectors $v\in \langle e_5, e_6 \rangle \cap \Lambda$ with $\| v \|^2=8\,$ are 
$\,v=\pm 2 e_5 \pm 2e_6\,$ and one has that:
$$\multline
\sum_{\text{ord}(B)=4} e^{2\pi i \, b \cdot (\pm 2 e_5 \pm 2 e_6)} = 
e^{2\pi i \, \frac{e_5}4 \cdot (\pm 2 e_5 \pm 2 e_6)} +
e^{2\pi i \, \frac{3e_5}4 \cdot (\pm 2 e_5 \pm 2 e_6)} + \\
+ e^{2\pi i \, (\frac {e_5}4+\frac{e_6}2 ) \cdot (\pm 2 e_5 \pm 2 e_6)} +
e^{2\pi i \, (\frac{3e_5}4+\frac{e_6}2 ) \cdot (\pm 2 e_5 \pm 2 e_6)} = -4.
\endmultline$$
Therefore  $\,\sum_{\text{ord}(B)=4} {{e}}_{8,B}(\Gamma) = -16$.
This  concludes the proof of the non $p$-isospectrality for $\,1\le p\le 5$.
\endexample

\example{Remark 5.2}
We note also that the non isospectrality in the previous example can be obtained by comparing $d_{p,0}(\Gamma)\,$ and $\,d_{p,0}(\Gamma')$, which can be determined directly by 
Remark 3.4.

One has that $d_{p,0}(\Gamma)=\beta_p(M_\Gamma)=\dim{\Lambda^p(\R^n)}^F$ and analogously for $d_{p,0}(\Gamma')$. 
The calculation of the $F$ and $F'$-invariants gives respectively: 

\

\centerline{
\vbox{\tabskip=0pt
\def\tablerule{\noalign{\hrule}}
\halign to280pt {\strut#& \vrule#\tabskip=1em plus 2em&
 \hfil# \hfil&  \vrule#& 
 \hfil# \hfil&  #& 
 \hfil# \hfil&  #& 
 \hfil# \hfil&  #& 
 \hfil# \hfil&  #& 
 \hfil# \hfil&  #& 
 \hfil# \hfil&  #& 
 \hfil# \hfil&  \vrule#  
\tabskip=0pt\cr\tablerule
&&  $p$\kern -10pt && $0$ && $1$ && $2$ && $3$  && $4$ && $5$ && $6$ & \cr\tablerule 
{\vrule height12pt width-6pt depth6pt} &&  $\beta_p(M_\Gamma)$ \kern -6pt && $1$ && $2$  && $3$ && $4$ && $3$ && $2$ && $1$  &\cr\tablerule
{\vrule height12pt width-6pt depth6pt} && $\beta_p(M_{\Gamma'})$ \kern -6pt&& $1$ && $1$ && $1$ && $2$ && $1$ && $1$ && $1$ &\cr
\tablerule
\cr}}
}
\endexample

\example{Example 5.3}
One notes that the group $\Gamma$ defined in Example 5.1 has the property that the holonomy action commutes with the
invariant K\"ahler structure   defined by 
$\,{\Cal J}_1=\ovr{\undr{\left[ \smallmatrix \widetilde J & & \\ & \widetilde J & \\ & & \widetilde J \endsmallmatrix \right]}}\,$ in $\R^6$.
It follows that $M_\Gamma$ inherits a K\"ahler structure. On the other hand $M_{\Gamma'}$ (as in 5.1) is isospectral to
$M_\Gamma$ but the holonomy action does not commute with ${\Cal J}_1$ and actually, $M_{\Gamma'}$  does not admit any
K\"ahler structure since $\beta_1(M_{\Gamma'})=1$ is odd (see \cite{We}).

One may use the previous pair to obtain a hyperk\"ahler manifold isospectral to a non-hyperk\"ahler one, as follows.
By duplication of the tables of ${\Gamma}$ and $\Gamma'$ in Example 5.1 (graphically, 
this means placing an identical second copy of each table below the first one) 
we obtain Bieberbach groups $\Gamma\!_2$ and $\Gamma'\!_2$ of dimension 12, with holonomy groups ${\Z}_4\times {\Z}_2$ and
holonomy representations $\rho \oplus \rho$ and $\rho' \oplus \rho'$, where $\rho$, $\rho'$ denote the holonomy representation
of $\Gamma$ and $\Gamma'$, respectively. By using the same bijection as before we see that $M_{\Gamma_2}$ and $M_{\Gamma'_2}$
are isospectral.

We extend the complex structure ${\Cal J}_{1}$  to $\R^{12}$ by
$\,{\Cal J}_{2,1}=\ovr{\undr{\left[ \smallmatrix {\Cal J}_{1} & 0\\0 & -{\Cal J}_{1}
\endsmallmatrix \right]}}\,$, and furthermore we have a second complex structure $\,{\Cal J}_{2,2}=\ovr{\undr{\left[
\smallmatrix 0& \I_6 \\  -\I_6 & 0 \endsmallmatrix \right]}}\,$ on $\R^{12}$, which anticommutes with 
${\Cal J}_{2,1}$, hence ${\Cal J}_{2,1}$ and  ${\Cal J}_{2,2}$ define a hyperk\"ahler structure on $\R^{12}$.
Since the holonomy action of $\Gamma_2$ clearly commutes with ${\Cal J}_{2,1}$ and  ${\Cal J}_{2,2}$,  $M_{\Gamma_2}$ now inherits a hyperk\"ahler structure. On the other hand, $\beta_1(M_{\Gamma'_2})= 2$,    the dimension of the fixed space of the holonomy action. 
Hence $M_{\Gamma'_2}$ can not carry a hyperk\"ahler structure, since $\beta_1$ is not divisible by~4.

\endexample

\example{Remark 5.4}
We may extend the holonomy representations $\,\rho$ and $\rho'$ of $F$ and $F'$ respectively, in Example 5.1   to
$\,\rho\oplus\tau$ and $\rho'\oplus \tau$ res\-pec\-ti\-ve\-ly, where $\tau$ is a sum of characters $\chi_h,\,\, 1\le
h\le k$, with values in $\{1,-1\}$. If we keep the same $b_i ,\, b'_i ,\,\, 1\le i\le 2$, clearly the resulting groups
$\Gamma$ and $\Gamma'$ are  torsion-free and of dimension $k+6$. Visually what has been done is adding $k$ rows of 1's
and $-1$'s  to the table in Example 5.1 (corresponding to the characters $\chi_1,\dots,\chi_k$). If, furthermore, one
chooses  $\chi_h$ in such a way that the  action of the second generator is trivial (that is, in such a way that all the
new entries in the fourth column are 1's), then $M_\Gamma$ and $M_{\Gamma'}$ are isospectral by  Theorem 3.1,  with the
bijection $\Phi$ chosen as in  Example 5.1.
\endexample

\example{Example 5.5}
If, as explained in Remark 5.4 one adds to the table in Example 5.1 only one character $\chi$ represented by the row
$(-1,1,-1,1,-1,1,-1)$ one obtains isospectral manifolds  $M_\Gamma$ and $M_{\Gamma'}$ of dimension 7, both non
orientable, with holonomy group $\Z_4\times\Z_2$, and which are not $p$-isospectral for $p=1,2,5,6$. This follows from
the values of the Betti numbers, which are given by: 

\

\centerline{
\vbox{\tabskip=0pt
\def\tablerule{\noalign{\hrule}}
\halign to320pt {\strut#& \vrule#\tabskip=1em plus 2em&
 \hfil# \hfil&  \vrule#& 
 \hfil# \hfil&  #& 
 \hfil# \hfil&  #& 
 \hfil# \hfil&  #& 
 \hfil# \hfil&  #& 
 \hfil# \hfil&  #& 
 \hfil# \hfil&  #& 
 \hfil# \hfil&  #& 
\hfil# \hfil&  \vrule#  
\tabskip=0pt\cr\tablerule
&&  $p$\kern -0pt && $0$ && $1$ && $2$ && $3$  && $4$ && $5$ && $6$ && $7$ & \cr\tablerule 
{\vrule height12pt width-6pt depth6pt} 
&&  $\beta_p(M_\Gamma)$ \kern -6pt && $1$ && $2$  && $3$ && $4$ && $3$ && $2$ && $1$ && $0$ &\cr\tablerule
{\vrule height12pt width-6pt depth6pt} 
&& $\beta_p(M_{\Gamma'})$ \kern -6pt&& $1$ && $1$ && $2$ && $4$ && $3$ && $3$ && $2$ && $0$ &\cr
\tablerule
\cr}}
}

It is easy  to verify that $M_\Gamma$ and $M_{\Gamma'}$ are $p$-isospectral for $p=3$, by the first statement in Remark
3.2, hence also for
$p=4$ and $7$, by the first statement in Remark 3.3. We note that both manifolds satisfy $\beta_p(M)\ne \beta_{n-p}(M)$,
hence spec$^p(M)\ne$ spec$^{n-p}(M)$, for every $p$, in contrast to the orientable case (see Remark 3.4).

\endexample

\example{Example 5.6}
In the present case we add  the character  $(-1,1,-1,1,-1,1,-1)$ twice, to $\Gamma$ and
the characters  $(-1,1,-1,1,-1,1,-1)$ and $(1,1,1,1,1,1,1)$ to $\Gamma'$ in Example 5.1,  respectively. 
One thus gets the following tables:

\

\centerline{
\vbox{\tabskip=0pt
\def\tablerule{\noalign{\hrule}}
\halign to340pt {\strut#& \vrule#\tabskip=1em plus 2em&
 \hfil# \hfil&  #& 
 \hfil# \hfil&  #& 
 \hfil# \hfil& \vrule#& 
 \hfil# \hfil& \vrule#& 
 \hfil# \hfil&  #&  
 \hfil# \hfil&  #& 
 \hfil# \hfil&  \vrule#  
\tabskip=0pt\cr\tablerule
&& {\vrule height14pt width0pt depth8pt} $B_1$ && $B_1^2$ && $B_1^3$ && $B_2$  && $B_1B_2$ && $B_1^2B_2$ &&  $B_1^3B_2$ & \cr\tablerule 
&& {\vrule height12pt width0pt depth0pt}  $\widetilde J$  && $-\I_2$  && $-\widetilde J$ && $-\I_2$ && $-\widetilde J$ && $\I_2$ && $\widetilde J$ &\cr
&&  $\widetilde J$  && $-\I_2$  && $-\widetilde J$ && $\I_2$ && $\widetilde J$ && $-\I_2$ && $-\widetilde J$ &\cr
{\vrule height0pt width0pt depth8pt}&& $1_{\frac 14}$ && $1_{\frac 12}$ && $1_{\frac 34}$ && 1 && 
$1_{\frac 14}$ && $1_{\frac 12}$ && $1_{\frac 34}$ &\cr
{\vrule height0pt width0pt depth8pt}&& 1 && 1 && 1 && $1_{\frac 12}$ && $1_{\frac 12}$ && $1_{\frac 12}$ && $1_{\frac 12}$ &\cr\tablerule
&& $-1$ && 1  && $-1$ && $1$ && $-1$ && $1$ && $-1$ &\cr
&& $-1$ && 1  && $-1$ && $1$ && $-1$ && $1$ && $-1$ &\cr
\tablerule
\cr}}
}

\centerline{
\vbox{\tabskip=0pt
\def\tablerule{\noalign{\hrule}}
\halign to348pt {\strut#& \vrule#\tabskip=1em plus 2em&
 \hfil# \hfil&  #& 
 \hfil# \hfil&  #& 
 \hfil# \hfil& \vrule#& 
 \hfil# \hfil& \vrule#& 
 \hfil# \hfil&  #&  
 \hfil# \hfil&  #& 
 \hfil# \hfil&  \vrule#  
\tabskip=0pt\cr\tablerule
&& {\vrule height14pt width0pt depth8pt} $B'_1$ && ${B'_1}^2$ && ${B'_1}^3$ && $B'_2$  && $B'_1B'_2$ && ${B'_1}^2B'_2$ &&  ${B'_1}^3B'_2$ & \cr\tablerule 
&& {\vrule height12pt width0pt depth0pt}  $\widetilde J$  && $-\I_2$  && $-\widetilde J$ && $-\I_2$ && $-\widetilde J$ && $\I_2$ && $\widetilde J$ &\cr
&& 1 && 1  && 1 && $-1$ && $-1$ && $-1$ && $-1$ &\cr
{\vrule height0pt width0pt depth8pt}&& $-1$ && 1  && $-1$ && $1_{\frac 12}$ && $-1_{\frac 12}$ && $1_{\frac 12}$ && $-1_{\frac 12}$ &\cr
{\vrule height0pt width0pt depth8pt}&& $-1$ && $1$ && $-1$ && 
$-1_{\frac 12} $ && $1_{\frac 12}$ && $-1_{\frac 12}$ && 
$1_{\frac 12}$ &\cr
{\vrule height0pt width0pt depth8pt}&& $1_{\frac 14}$ && $1_{\frac 12}$ && $1_{\frac 34}$ && $1$ && $1_{\frac 14}$ && $1_{\frac 12}$ && $1_{\frac 34}$ &\cr\tablerule
&& $-1$ && 1  && $-1$ && $1$ && $-1$ && $1$ && $-1$ &\cr
&& 1 && 1  && 1 && $1$ && $1$ && $1$ && $1$ &\cr
\tablerule
\cr}}
}

We see that $M_\Gamma$ is orientable while $M_{\Gamma'}$ is not. We shall now show that $M_\Gamma$ and $M_{\Gamma'}$ are isospectral on $p$-forms if and only if $p$ is odd.

We shall use an argument similar to that of Example 5.1 with the same bijection 
$\Phi : F \rightarrow F'$. We have again that, for each $p$, the contribution of $B_1^2$  to $d_{p,\mu}(\Gamma)$ in (3.1)
is the same as the contribution of $B'_2$  to $d_{p,\mu}(\Gamma')$ and in the same way, the contributions of $B_2$ and
${B'_1}^2$ are the same, and also those of $B_1^2B_2$ and ${B'_1}^2B'_2$. Therefore, the analysis of  isospectrality is
reduced to the comparison of the contributions of the elements of order 4 in $F$ and $F'$. We also observe that tr$_p(B)$
(resp\. tr$_p(B')$) is the same for all 4 elements of order 4 in $F$ (resp\. $F'$), for each $p$. Hence, $M_\Gamma$ and
$M_{\Gamma'}$ will be $p$-isospectral if and only if, for each $\mu$,
 
$$\text{tr}_p(B_1)\sum_{B:\text{ord}(B)=4} {{e}}_{\mu,B}(\Gamma)=
\text{tr}_p(B'_1)\sum_{B':\text{ord}(B')=4} {{e}}_{\mu,B'}(\Gamma').\tag5.1$$

We get the following  values of the $p$-traces, whose verification will be  omitted since it is similar to the  calculation in  Example
5.1:

\

\centerline{
\vbox{\tabskip=0pt
\def\tablerule{\noalign{\hrule}}
\halign to320pt {\strut#& \vrule#\tabskip=1em plus 2em&
 \hfil# \hfil&  \vrule#& 
 \hfil# \hfil& \kern-7pt  #& 
 \hfil# \hfil& \kern-7pt  #& 
 \hfil# \hfil& \kern-7pt  #& 
 \hfil# \hfil& \kern-7pt  #& 
 \hfil# \hfil& \kern-7pt  #& 
 \hfil# \hfil& \kern-7pt  #& 
 \hfil# \hfil& \kern-7pt  #& 
 \hfil# \hfil& \kern-7pt  #& 
 \hfil# \hfil& \kern-7pt  \vrule#  
\tabskip=0pt\cr\tablerule
&&  $p$\kern -0pt && 0&& $1$ && $2$ && $3$  && $4$ && $5$ && 6 && 7 && 8 & \cr\tablerule 
&&  tr$_p(B_1)$ \kern -6pt && 1 && 0 && 0 && 0 && $-2$  && $0$ && $0$ && $0$ && $1$  &\cr\tablerule
&& tr$_p(B'_1)$ \kern -6pt&& 1 && $0$ && $-2$ && $0$ && 0 && 0 && $2$ && $0$ && $-1$ &\cr
\tablerule
\cr}}
}

Since tr$_p(B_1)=\text{tr}_p(B'_1)=0,\,\,\,\forall \, p$ odd, it turns out  that $M_\Gamma$ and $M_{\Gamma'}$ are isospectral on $p$-forms for all $p$ odd. We now show that they are not $p$-isospectral for $p$ even.

$$\multline
\sum_{\Sb {j=5,6} \\ {B : \text{ord}(B)=4}\endSb} e^{2\pi i \, b \cdot (\pm e_j)} = \sum_{j=5,6} \bigg( e^{2\pi i \, \frac{e_5}4 \cdot (\pm e_j)} + 
e^{2\pi i \, \frac{3e_5}4 \cdot (\pm e_j)} + \\
+ e^{2\pi i \, (\frac {e_5}4+\frac{e_6}2 ) \cdot (\pm e_j)} +
e^{2\pi i \, (\frac{3e_5}4+\frac{e_6}2) \cdot (\pm e_j)} \bigg)= 0. 
\endmultline$$
$$\sum_{\Sb {B' : \text{ord}(B')=4} \\ j:B'e_j=e_j \endSb} 
e^{2\pi i \, b' \cdot (\pm e_j)} = 
\sum_{\Sb B' : \text{ord}(B')=4 \\ j:B'e_j=e_j \\ 1\le j\le 6\endSb} e^{ 2\pi i \, b' \cdot (\pm e_j)} +
\sum_{B' : \text{ord}(B')=4} e^{ 2\pi i \, b' \cdot (\pm e_8)} = 0+4=4.
$$

Therefore one concludes that $\dsize\undr{\sum_{B : \text{ord}(B)=4} {{e}}_{1,B}(\Gamma)=0}\,\,\,$ and $\,\dsize\undr{\sum_{B' : \text{ord}(B')=4} {{e}}_{1,B'}(\Gamma')=8}\,\,$ and according to the tables giving tr$_p(B_1)$ and tr$_p(B'_1)$ one obtains 
that $M_\Gamma$ and $M_{\Gamma'}$ are not $p$-isospectral for $p=0, 2, 6$ and 8. For $p=4$ one can take $\mu=8$ and using the calculations in  Example
5.1, one obtains
$$\sum_{\text{ord}(B)=4} {{e}}_{8,B}(\Gamma) = -16.$$ 
Since tr$_4(B_1)=-2\,$ and $\,\text{tr}_4(B'_1)=0\,$ it turns out that $\,d_{4,8}(\Gamma)\ne d_{4,8}(\Gamma')\,$ hence $\,M_\Gamma$ and $M_{\Gamma'}\,$ are not isospectral on $4$-forms.

Again, one could also have concluded  the non $p$-isospectrality of $M_\Gamma$ and $M_{\Gamma'}$ for $p$ even, $p>0$ by
computing their Betti numbers:

\

\centerline{
\vbox{\tabskip=0pt
\def\tablerule{\noalign{\hrule}}
\halign to280pt {\strut#& \vrule#\tabskip=1em plus 2em&
 \hfil#  \hfil&  \vrule#& 
 \hfil# \hfil& \kern-7pt  #& 
 \hfil# \hfil& \kern-7pt  #& 
 \hfil# \hfil& \kern-7pt  #& 
 \hfil# \hfil& \kern-7pt  #& 
 \hfil# \hfil& \kern-7pt  #& 
 \hfil# \hfil& \kern-7pt  #& 
 \hfil# \hfil& \kern-7pt  #& 
 \hfil# \hfil& \kern-7pt  #& 
 \hfil# \hfil& \kern-7pt  \vrule#  
\tabskip=0pt\cr\tablerule
&&  $p$\kern -10pt && $0$ && $1$ && $2$ && $3$  && $4$ && $5$ && $6$ && 7 && 8 & \cr\tablerule 
{\vrule height12pt width-6pt depth6pt} &&  $\beta_p(M_\Gamma)$ \kern -6pt && $1$ && $2$  && $4$ && $6$ && $6$ && $6$ && $4$ && 2 && 1 &\cr\tablerule
{\vrule height12pt width-6pt depth6pt} && $\beta_p(M_{\Gamma'})$ \kern -6pt&& $1$ && $2$ && $3$ && $6$ && $7$ && $6$ && $5$ && 2 && 0 &\cr
\tablerule
\cr}}
}

\endexample

\example{Example 5.7}
If one takes $\Gamma$ as in  Example 5.6 and $\Gamma'$ by adding twice the trivial character to $\Gamma'$ in Example 5.1,
it turns out  that $M_\Gamma$ and $M_{\Gamma'}$ are $p$-isospectral for $p=2$ and 6 but not for the remaining  values of
$p$. Indeed, by arguing as in the previous examples, we see that $p$-isospectrality occurs if and only if (5.1) holds in
this case. Since one verifies that tr$_2(B_1)=\text{tr}_2(B'_1)=0$, then $M_\Gamma$ and $M_{\Gamma'}$ are isospectral on
2-forms. Since $M_\Gamma$ and $M_{\Gamma'}$ are both orientable, by Remark 3.3, they are isospectral on 6-forms. The non
isospectrality for $p\ne 2, 6\,$ can be verified as in the previous examples.
\endexample

\spb

\example{Example 5.8}
In the present example, we shall give two flat manifolds of dimension 4 which are isospectral on $p$-forms for $p$ odd and  having different holonomy groups: $\Z_2^2$ and $\Z_4$, respectively.  

Let $\Gamma = \langle B_1L_{b_1}, B_2L_{b_2} \Lambda \rangle$, and 
$\Gamma' = \langle B'L_{b'}, \Lambda \rangle$, where 
$\Lambda = \Z^4$ and 
$$
B_1=\left[ \matrix  1 & & & \\  & 1 & &\\ & & -1 & \\
& & & -1 \endmatrix \right] \qquad
B_2=\left[ \matrix  1 & & & \\  & -1 & &\\ & & -1 & \\
& & & 1 \endmatrix \right]\qquad \qquad
B'=\left[ \matrix \ovr {\undr{\widetilde J }}& &  \\  & -1 & \\ & & 
\ovr {\undr{ 1}} 
\endmatrix \right]
$$
$$b_1=\frac{e_1}2, \,\, b_2=\frac{e_4}2,\qquad b'=\frac{e_4}4.$$ 
By a verification of the conditions (i) and (ii) in Proposition 2.1,
as in Example 5.1 (but simpler), we get that
$\Gamma$ and $\Gamma'$ are Bieberbach groups with holonomy groups $\Z_2 \times \Z_2$ and $\Z_4$, respectively.

We shall give two tables listing the non trivial elements in $F$ and $F'$, together with subindices indicating the non zero translational components: for instance, since $b_1=\frac{e_1}2$, we write $\frac{1}2$ as a subindex of the first diagonal element of $B_1$.

\

\centerline{
\vbox{\tabskip=0pt
\def\tablerule{\noalign{\hrule}}
\halign to136pt {\strut#& \vrule#\tabskip=1em plus 2em&
 \hfil# \hfil&  #& 
 \hfil# \hfil&  #& 
 \hfil# \hfil& \vrule# 
 \tabskip=0pt\cr\tablerule
&& {\vrule height14pt width0pt depth8pt} $B_1$ && $B_2$ && $B_1B_2$  & \cr\tablerule 
&& {\vrule height12pt width0pt depth0pt}  $1_{\frac 12}$  && $1$  && $1_{\frac 12}$  &\cr
&&  $1$  && $-1$  && $-1$ &\cr
&&  $-1$  && $-1$  && $1$ &\cr
&&  $-1\undr{}$  && $1_{\frac 12}$  && $-1_{\frac 12}$ &\cr
\tablerule
\cr}}
\qquad
\vbox{\tabskip=0pt
\def\tablerule{\noalign{\hrule}}
\halign to140pt {\strut#& \vrule#\tabskip=1em plus 2em&
 \hfil# \hfil&  #& 
 \hfil# \hfil&  #& 
 \hfil# \hfil& \vrule# 
 \tabskip=0pt\cr\tablerule
&& {\vrule height14pt width0pt depth8pt} $B'$ && ${B'}^2$ && ${B'}^3$  & \cr\tablerule 
&& {\vrule height12pt width0pt depth0pt}  $\ovr{\undr{\widetilde J}}$  && $-\I_2$  && $-\widetilde J$  &\cr
&&  $\ovr{\undr{-1}}$  && $1$  && $-1$ &\cr
{\vrule height0pt width0pt depth8pt}&& $\ovr{1_{\frac 14}}$ && $1_{\frac 12}$ && $1_{\frac 34}$  &\cr
\tablerule
\cr}}
}

It is not hard to verify that $\text{tr}_p (B)=0$ for any $B\in F, F'$, $B\ne \I$ and $p=1,3$, hence Theorem 3.1  implies
that $M_\Gamma$ and $M_{\Gamma'}$ are $p$-isospectral for $p=1,3$. 

By using the  above pair as a starting point, we have constructed, for each $m\ge 2$, a pair of flat $2m$-manifolds with holonomy groups $\Z_2^m$  and $\Z_4\times\Z_2^{m-2}$ respectively, which are $p$-isospectral for any $p$ odd. For brevity we shall omit their description.

\endexample

\

\example{Example 5.9}
Adding, as in Remark 5.4, $k$ trivial characters to the groups in Example 5.1, one obtains 0-isospectral manifolds
$\,M_\Gamma\times T^k\,$ and $\,M_{\Gamma'}\times T^k$, of dimension $n=k+6$. Calculating the Betti numbers by the
K\"unneth formula  one gets:
$$\beta_h(M_\Gamma\times T^k)=\sum_{i=0}^{\max (h,6)} \beta_i(M_\Gamma) \binom k{h-i},$$
and similarly for $M_{\Gamma'}\times T^k$, under the convention that $\,\binom nm=0$, if $\,m>n$.
{}From the  table of Betti numbers in Example 5.1 one concludes that 
$$\beta_h(M_\Gamma\times T^k) > \beta_h(M_{\Gamma'}\times T^k)\tag5.2$$
for $0<h<n$. As a consequence, for any $\,n\ge 6$, one obtains  manifolds which are isospectral on functions but not on $p$-forms for $\,0<p<n$.
\endexample

\example{Remark 5.10}
As it was  already mentioned, examples of  isospectral, non homeomorphic manifolds, which 
are not isospectral on $p$-forms have been constructed by A.\ Ikeda (see \cite{Ik})  for lens spaces ($n\ge 5$), by C.\
Gordon, R.\ Gornet (see \cite{Go}, \cite{Gt}) in the context of nilmanifolds ($n\ge 5$) and D.\ Schueth  (\cite{Sch}). 

It is an open question  whether, for every given subset $I$ of $\{0,1,\dots,n\}$, there exist compact manifolds which are $p$-isospectral if and only if $p \in I$.  For each $k\ge 0$, Ikeda has constructed lens spaces which are 
isospectral on $p$-forms for $0\le p\le k$, but not  on $(k+1)$-forms.
Also, we have seen that flat manifolds allow to give many examples of pairs which are isospectral for certain values of $p$ only. However, even though more examples can be given,  the construction of  flat manifolds which are $p$-isospectral for only a fixed set of values of $p$, is still a complicated matter, since one has to keep control of all values of the $p$-traces, $\text{tr}\!_p$, at the same time.

We have also seen  that the  isospectral manifolds given in many of the  examples above are topologically quite different, in particular for those in Example
5.9 one has $\beta_p(M_\Gamma) > \beta_p(M_{\Gamma'})\, \text{ for } \, 0 < p < n$ (hence $\text{spec}^p\,M_{\Gamma} \ne
\text{spec}^p\, M_{\Gamma'}$, for $0<p<n$). We note that all spherical space forms (hence all Ikeda's lens spaces) have
the same de Rham cohomology as the sphere $S^n$ and their topological distinction is more delicate (see \cite{Gi}). Also,
almost all of the isospectral  nilmanifolds in the references mentioned above are of the form
$\Gamma\backslash N$, $\Gamma'\backslash N$, where $\Gamma$ and $\Gamma'$ are lattices in the same simply connected nilpotent Lie group $N$, hence, by a theorem of Nomizu, their real cohomology coincides with the Lie algebra cohomology of ${\goth n}$, the Lie algebra of $N$ (hence  they both have the same Betti numbers). To our best knowledge, examples of $0$-isospectral nilmanifolds with Betti numbers showing property
(5.2) (hence not $p$-isospectral for  $0<p<n$) have not yet been given. Their construction would involve finding an
isospectral pair  $\Gamma\backslash N$, $\Gamma'\backslash N'$ with $N\ne N'$, where one can compute, or at least
compare,  the dimensions of  $H^p({\goth n})$  and $H^p({\goth n}')$ for each value of $p$ and this would not seem such a
simple matter. 

\endexample

\example{Acknowledgments}
The authors wish to thank Carolyn Gordon and Peter Gilkey for very useful comments on an earlier 
version of this article and Alex Samorodnitsky for useful conversations on the Krawtchouk polynomials.
\endexample

\enddocument